\documentclass[a4,11pt]{article}
\usepackage{calc}
\usepackage[all]{xy}
\usepackage[centertags]{amsmath}
\usepackage{latexsym}
\usepackage{amsfonts}
\usepackage{amssymb}
\usepackage{amsthm}
\usepackage[all]{xy}
\usepackage{fancyhdr}
\usepackage [dvips]{epsfig}

\usepackage{newlfont}
\usepackage[latin1]{inputenc}
\usepackage[french,english]{babel}
\usepackage{graphicx}
\usepackage{t1enc}
\theoremstyle{plain}

\newtheorem{propo}{Proposition}[section]
\newtheorem{thm}{Theorem}[section]
\newtheorem{defi}{Definition}[section]
\newtheorem{coro}{Corollary}[section]

\theoremstyle{plain}
\begin{document}
\title{On the Hochschild homology of elliptic Sklyanin algebras}
\author{S. R. Tagne Pelap}
\date{May 19, 2008}
\maketitle
\begin{abstract}
In this paper, we compute the Hochschild homology of elliptic Sklyanin algebras. These algebras are deformations of polynomial algebra with a Poisson bracket called the Sklyanin Poisson bracket.\\
\end{abstract}
\nocite{*}
\section*{Introduction}
The family of algebras defined by Sklyanin in (~\cite{sky1}), which today carries his name, is naturally associated with two parameters: an elliptic curve and a point on this curve. It is a family of associative algebras with $4$ generators and six quadratic relations. These algebras are flat deformations of polynomial algebras with four variables.\\
 The paper (~\cite{ode1}) of  Odesskii and Feigin gives a generalization of these algebras. Considering an elliptic curve $\mathcal E$ and a point $\eta$ of this curve, the elliptic algebras or algebras of Feigin-Odesskii are the family of associative algebras $Q_{n, k}(\mathcal E, \eta)$, $k<n$ (which are mutually prime), with $n$ generators and the relations:
 \begin{equation}
 \displaystyle\sum_{r\in\mathbb{Z}/n\mathbb{Z}}\frac{\theta_{j-i+r(k-1)}(0)}{\theta_{kr}(\eta)\theta_{j-i-r}(-\eta)}x_{j-r}x_{i+r}=0
 \end{equation}
 where $\theta_{\alpha}$, $\alpha\in\mathbb{Z}/n\mathbb{Z}$, are theta functions (~\cite{ode1}).\\
  These algebras arise as deformations of Poisson algebras.\\
 The goal of this paper is to compute the Hochschild homology of $Q_{4, 1}(\mathcal E, \eta).$\\
 In (\cite{kon}), Kontsevich proved that there is an isomorphism between Poisson cohomology of a Poisson algebra and Hochschild cohomology of a particular quantization of this Poisson algebra. This result is also known as Kontsevich's formality theorem. Covariant versions of this formality theorem were conjectured by B. Tsygan (\cite{tsy}) also for this particular quantization and proved by B. Shoikhet (\cite{sho}). Since it is unknown if our quantization is the same as the one given by the formality theorem, we can't use this result. To compute our homology, we choose follow a similar method as Van den Bergh in his computation of Hochschild homology of the algebra $Q_{3, 1}(\mathcal E, \eta).$\\
 The paper is organized as follows. We start by reviewing some general facts on the Hochschild homology, Poisson homology and, Koszul algebras.
The next part is devoted to generalities on associative elliptic algebra $Q_{4, 1}(\mathcal E, \eta),$ also called the Sklyanin algebra. The last and the main part of this paper is devoted to the computation of the Hochschild homology of the Sklyanin algebra $Q_{4, 1}(\mathcal E, \eta).$ \\
 {\bf Acknowledgements.} The author is grateful to Jean-Claude Thomas for useful comments and discussions. This work is a part of my thesis prepared at the University of Angers. I would like to take this opportunity to thank my advisors, Vladimir Roubtsov and Bitjong Ndombol, for suggesting that I tackle this interesting problem and I also thank them for being available during this projet. I thank the referee of this paper for his comments and suggestions which are incorporated in this revised version. This work was partially supported by the Programme SARIMA.
\section{Preliminary facts}
\subsection{Hochschild homology}
Let $A$ be an associative algebra over $K$, where $K$ is a field.\\
One can define the complex $(C_\star(A), b)$ by:
$$C_n(A)=A^{\otimes(n+1)};$$
$$b(a_0\otimes\cdot\cdot\cdot\otimes a_n)=(-1)^na_na_0\otimes\cdot\cdot\cdot\otimes a_{n-1}+\displaystyle{\sum_{i=0}^{n-1}}(-1)^ia_0\otimes\cdot\cdot\cdot\otimes a_ia_{i+1}\otimes\cdot\cdot\cdot\otimes a_n$$
The Hochschild homology of $A$ with coefficients in $A$ is the homology of the complex $(C_\star(A), b).$ This homology is denoted by $HH_\star(A).$\\
We denote by $B$ Connes's coboundary $C(A)\longrightarrow C(A)$ defined as follows:
$$B(a_0\otimes\cdot\cdot\cdot\otimes a_n)=\displaystyle{\sum_{i=0}^{n}}(-1)^{ni}1\otimes a_i\otimes\cdots\otimes a_{n}\otimes a_0\otimes\cdots\otimes a_{i-1}+(-1)^n\displaystyle{\sum_{i=0}^{n}}(-1)^{ni}a_i\otimes\cdots\otimes a_{n}\otimes a_0\otimes\cdots\otimes a_{i-1}\otimes 1.$$
We have $b\circ B +B\circ b=0.$
\subsection{Poisson homology}
A Poisson bracket on a commutative algebra $\mathcal R$ is an antisymmetric biderivation $\{\cdot, \cdot\}: \mathcal R\times\mathcal R\rightarrow\mathcal R$ such that $(\mathcal R, \{\cdot, \cdot\})$ is a Lie algebra.Then $(\mathcal R, \{\cdot, \cdot\})$ is called a Poisson algebra. One can also say that  $\mathcal R$ is endowed with a Poisson structure.The elements of the center of a Poisson algebra $(\mathcal R, \{\cdot, \cdot\})$ (that means the elements $a\in\mathcal R$ such that $\{a, b\}=0$ for all $b\in\mathcal R$) are called the Casimirs of this Poisson algebra.\\
Let us give an example: consider $n-2$ polynomials $P_i$ in $K[x_1,\cdots, x_n],$ where $K$ is a field of characteristic zero. For any polynomial $\lambda\in K[x_1,...,x_n]$, we can define a bilinear application:\\
$$\{\cdot ,\cdot\} : K[x_1,...,x_n]\otimes K[x_1,...,x_n]\longrightarrow K[x_1,...,x_n]$$
by the formula
\begin{equation}\label{q}
\{f,g\}=\lambda\frac{df\wedge dg\wedge dP_1\wedge...\wedge dP_{n-2}}{dx_1\wedge dx_2\wedge...\wedge dx_n},\ \ f,g\in K[x_1,...,x_n]
 \end{equation}
This operation gives a Poisson algebra structure on $K[x_1,...,x_n],$ also called a Jacobian Poisson structure (JPS). Naturally the polynomials $P_i, i=1,...,n-2$ are Casimir functions for the brackets (\ref{q}) and every Poisson structure in $K^n$ with $n-2$ generic Casimirs $P_i$ is written in this form. \\
The case $n=4$ in $(\ref{q})$ corresponds to the classical (generalized) Sklyanin Poisson algebra. The real Sklyanin algebra is associated with the following two quadrics in $K^4$:
\begin{equation}\label{e}
P_1=x_1^2+x_2^2+x_3^2
\end{equation}
\begin{equation}\label{f}
P_2=x_4^2+J_1x_1^2+J_2x_2^2+J_3x_3^2
\end{equation}
\\
Let $\{\cdot, \cdot\} : \mathcal R\times\mathcal R\rightarrow\mathcal R$ be a Poisson bracket on an algebra $\mathcal R.$\\
The Poisson boundary operator, also called the Brylinski or Koszul differential and denoted by $$\partial : \Omega^\bullet(\mathcal R)\longrightarrow\Omega^{\bullet-1}(\mathcal R)$$ is given by:
$$\partial_k(F_0dF_1\wedge...\wedge dF_k)=\displaystyle{\sum_{1\leq i\leq k}}(-1)^{i+1}\{F_0, F_i\}dF_1\wedge...\wedge \widehat{dF_i}\wedge ...\wedge dF_k$$
$$+\displaystyle{\sum_{1\leq i<j\leq k}} (-1)^{i+j}F_0d\{F_i,F_j\}\wedge dF_1\wedge...\wedge \widehat{dF_i}\wedge...\widehat{dF_j}\wedge...\wedge dF_k$$\\
where $F_0,...,F_k\in\mathcal R.$\\
One can check, by a direct computation, that $\partial_k$ is well-defined and that it is a boundary operator.\\
The homology of the complex $(\Omega^\bullet(\mathcal R), \partial)$, denoted by $PH_\star(\mathcal R, \partial)$, is called the Poisson homology associated with the Poisson bracket $\{\cdot, \cdot\}.$ \\
In (~\cite{van}) Michel Van den Bergh computes the Poisson homology of $q_{3,1}(\mathcal E)$ which is the Jacobian Poisson structure given by the polynomial $P=\frac{1}{3}(x_1^3+x_2^3+x_3^3)+kx_1x_2x_3.$
Using the similar method as Van den Bergh, Nicolas Marconnet computes the Poisson homology of a cubic Jacobian Poisson structure on the polynomial algebra $K[x_1, x_2, x_3]$ given by a polynomial $\phi=\frac{1}{2}x_3^2+\frac{q_1}{4}x_1^4+\frac{q_1}{4}x_2^4+\frac{3q_1}{2}x_1^2x_2^2.$
The another point of view have been considered by Anne Pichereau in (~\cite{pic}): she computes the Poisson homology of Jacobian Poisson structures in dimension three given by a weight homogeneous polynomial with an isolated singularity.\\
In our article (~\cite{pelap1}), we compute the Poisson homology of a Jacobian Poisson structure in dimension four given by weight homogeneous polynomials $P_1$ and $P_2$ which form a complete intersection. We proved that the Poincaré series of these homological groups depend only on the weights and the degrees of $P_1$ and $P_2.$  In the quadratic homogeneous case, we obtained the following result:
\begin{propo}\textnormal{(~\cite{pelap1})}
The Poincaré series of the Poisson homological groups of Jacobian Poisson structures in dimension four given by quadratic homogeneous polynomials $P_1$ and $P_2$ which form a complete intersection, $PH_i(\mathcal R), \partial), i=0, 1, 2, 3, 4,$ have the following forms:\\
$$\begin{array}{lcl}
P(PH_0(\mathcal R, \partial),t)&=&\frac{2t^2+4t+1}{(1-t^2)^2} ;\\
&&\\
P(PH_1(\mathcal R, \partial),t)&=&\frac{t^4+4t^3+4t^2+4t}{(1-t^2)^2};\\
&&\\
P(PH_2(\mathcal R, \partial),t)&=&\frac{2t^4+4t^3}{(1-t^2)^2}; \\
&&\\
P(PH_3(\mathcal R, \partial),t)&=&\frac{t^4}{(1-t^2)^2} ;\\
&&\\
P(PH_4(\mathcal R, \partial),t)&=&\frac{t^4}{(1-t^2)^2}. \\
\end{array}$$
\end{propo}
In Sklyanin's case, we obtained the following explicit result:
\begin{propo}\textnormal{(~\cite{pelap1})}\label{sph} In the generic case, the Poisson homology of the Sklyanin Poisson structure is described as follow:
\begin{enumerate}
             \item The homological group $PH_0(\mathcal R, \partial)$ is a rank 7 free $K[P_1, P_2]$-module generated by $(\mu_i)_{0\leq i\leq 6} = (1, x_1, x_2, x_3, x_4, x_1^2,x_3^2);$
             \item $PH_1(\mathcal R, \partial)$ is a rank 13 free $K[P_1, P_2]$ module given by:
$$PH_1(\mathcal R, \partial) \cong (\displaystyle{\bigoplus_{k=1}^6} K[P_1, P_2]d\mu_k)\oplus(\displaystyle{\bigoplus_{k=1}^5} K[P_1, P_2]\mu_kdP_1)\oplus K[P_1, P_2]dP_1\oplus K[P_1, P_2]dP_2;$$
             \item $PH_2(\mathcal R, \partial)$ is a rank 6 free $K[P_1, P_2]$ module given by:
$$PH_2(\mathcal R, \partial) \cong \left(\displaystyle{\bigoplus_{k=1}^5} K[P_1, P_2](d\mu_k\wedge dP_1)\right)\oplus K[P_1, P_2](dP_1\wedge dP_2);$$
             \item $PH_3(\mathcal R, \partial)$ is a rank $1$ free $K[P_1, P_2]$-module generated by $\pi$;
             \item $PH_4(\mathcal R, \partial)$ is a rank $1$ free $K[P_1, P_2]$-module generated by $\delta$,
\end{enumerate}
where             $\delta=dx_1\wedge dx_2\wedge dx_3\wedge dx_4$, and $\pi=x_1dx_2\wedge dx_3\wedge dx_4+x_2dx_3\wedge dx_1\wedge dx_4+x_3dx_1\wedge dx_2\wedge dx_4+x_4dx_2\wedge dx_1\wedge dx_3.$\\
\end{propo}
\subsection{Generalities on Koszul algebras}
Let $V$ be a finite -dimensional $K$-vector space and let $T(V)$ be the tensor algebra of $V$ over $K.$ Consider a quadratic $K$-algebra $A=T(V)/(W)$, where $W\subset V\otimes V.$ Let $W^\ast$ be the dual space of $W$ and $W^\perp\subset V^\ast\otimes V^\ast$ be the orthogonal of $W.$ The dual algebra of $A$ is defined as $A^!:=T(V^\ast)/(W^\perp).$\\
Let $(x_i)_{i=0,\cdots {n-1}}$ be a basis of $V$ and $(\zeta_i)_{i=0,\cdots {n-1}}$ its dual basis. Consider $e=\displaystyle{\sum_{i=0}^{n-1}}x_i\otimes\zeta_i\in A\otimes A^!.$ We have $e^2=0$ (~\cite{van}). Let $K_m(A)=A\otimes (A^!_m)^\ast$ and $K_\star (A)=\oplus_{m\geq 0}K_m(A).$ \\
If $f\in(V^{*\otimes n})^*$ and $x\in V^*$, then the inner product of $f$ by $x$, $f\cdot x\in(V^{*\otimes (n-1)})^*$, is defined by:
$$
\begin{array}{llcl}
f\cdot x :&(V^{*\otimes (n-1)})&\longrightarrow& k\\
&v_1\otimes\cdots\otimes v_{n-1}&\longmapsto& f(x\otimes v_1\otimes\cdots\otimes v_{n-1}).
\end{array}
$$
Then the right multiplication by $e$ (multiplication of algebra on $A$ and inner product on $(A^!)^*$) induces a map $d : K_m(A)\longrightarrow K_{m-1}(A)$ and a differential $d : K_\star A\longrightarrow K_\star A.$\\
The complex\\
\begin{equation}\label{kc}
A=A\otimes k=A\otimes (A^!_0)^\ast\stackrel{\times e}{\longleftarrow}A\otimes (A^!_1)^\ast\stackrel{\times e}{\longleftarrow}A\otimes (A^!_2)^\ast\stackrel{\times e}{\longleftarrow}\cdots
\end{equation}\\
is called the Koszul complex of $A.$\\
The augmented Koszul complex of $A$ is the complex:
\begin{equation}\label{akc}
0\longleftarrow k\stackrel{\varepsilon}{\longleftarrow} A=A\otimes k=A\otimes (A^!_0)^\ast\stackrel{\times e}{\longleftarrow}A\otimes (A^!_1)^\ast\stackrel{\times e}{\longleftarrow}A\otimes (A^!_2)^\ast\stackrel{\times e}{\longleftarrow}\cdots
\end{equation}
where $\varepsilon$ is the canonical projection.
\begin{defi}
A quadratic algebra $A$ is said to be a Koszul algebra if the augmented Koszul complex (\ref{akc}) is exact.
\end{defi}
\begin{propo}\textnormal{(~\cite{van})}
Suppose that $A$ is a Koszul algebra. Then $HH_\star(A)=H_\star(K(A), b).$
\end{propo}
Hence when $A$ is a Koszul algebra, we have two complexes which enable us to compute the Hochschild homology of $A$: $(K(A), b)$ and $(C(A), b)$.
Let us give an explicit quasi-isomorphism between them.\\
Since $(A_m^!)^\ast=\cap_{i+j+2=m}V^{\otimes i}\otimes W\otimes V^{\otimes j}$ (~\cite{van}), we define the map $q : A\otimes(A_m^!)^\ast\longrightarrow A\otimes A^{\otimes m}$ as being the restriction of the natural inclusion $A\otimes V^{\otimes m}\hookrightarrow A\otimes A^{\otimes m}.$
\begin{propo}\textnormal{(~\cite{van})}
$q : K(A)\longrightarrow C(A)$ is a quasi-isomorphism.
\end{propo}
\section{Generalities on non commutative Sklyanin algebras}
Here, we follow the initial description of Sklyanin algebras by Sklyanin in (~\cite{sky2}).
Let $\tau\in \mathbb{C}$ so that $Im(\tau)>0.$\\
Consider the subgroup $\Gamma=\mathbb Z\oplus\mathbb Z\tau.$\\
Let $\theta_{00},\ \theta_{01},\ \theta_{10},\ \theta_{11},$ from $\mathbb C$ to $\mathbb C,$ be Jacobi's theta functions  associated with $\Gamma$, as described in (~\cite{smsf}).\\
These holomorphic functions satisfy the following properties: \\
$$\theta_{ab}(z+1)=(-1)^{a}\theta_{ab}(z);$$
$$\theta_{ab}(z+\tau)=e^{(-\pi i\tau-2\pi iz-\pi ib)}\theta_{ab}(z);$$
and the zeros of $\theta_{ab}$ are the points: $\frac{1}{2}(1-b)+(1+a)\tau+\Gamma.$\\
Fix $\eta\in\mathbb C$ such that $\eta$ is not of order $4$ in $\mathbb C/\Gamma.$ Let $(ab, ij, kl)$ be a cyclic permutation of $(00, 01, 10, 11).$\\
Setting:
$$\alpha_{ab}=(-1)^{a+b}\left[\frac{\theta_{11}(\eta)\theta_{ab}(\eta)}{\theta_{ij}(\eta)\theta_{kl(\eta)}}\right]^2.$$ Consider $V$ a four dimensional $\mathbb C$-vector space with $S_0, S_1, S_2, S_3$ as a basis.\\
We define $A=T(V)/(I_2)$ where $(I_2)$ is the two-sided ideal generated by the subspace $I_2\subseteq V\otimes V$ with basis:\\
$$f_{0i}=[S_0, S_i]-\alpha_i(S_jS_k+S_kS_j);$$
$$f_{jk}=[S_j, S_k]-(S_0S_i+S_iS_0),$$
$\alpha_1=\alpha_{00},$ $\alpha_2=\alpha_{01},$ $\alpha_3=\alpha_{10},$ and $(i,j,k)$ is a cyclic permutation of $(1,2,3).$\\
A simple computation shows that $\alpha_1+\alpha_2+\alpha_3+\alpha_1\alpha_2\alpha_3=0.$\\
$A$ is called the Sklyanin algebra.\\
We have the following results from the paper of S.P. Smith and J.T. Stafford (~\cite{smsf}):
\begin{propo}\textnormal{(~\cite{smsf})}
$A$ is a Koszul algebra.
\end{propo}
Therefore the Hochschild homology of $A$ is given by the complex $(K(A), b)$, where $b$ is the Hochschid boundary.
\begin{propo}\textnormal{(~\cite{smsf})}
For $m\geq 5$, $A^!_m=0$ and if $m\leq4$, $A^!_m$ is spanned by the following elements:\\
 $A^!_0: 1$\\
 $A^!_1: \zeta_0, \zeta_1, \zeta_2, \zeta_3$\\
  $A^!_2: \zeta_0\zeta_1, \zeta_0\zeta_2, \zeta_0\zeta_3, \zeta_1\zeta_0, \zeta_2\zeta_0$\\
   $A^!_3: \zeta_0\zeta_1\zeta_0, \zeta_0\zeta_2\zeta_0, \zeta_0\zeta_3\zeta_0, \zeta_1\zeta_0\zeta_1$\\
   $A^!_4: \zeta_0\zeta_1\zeta_0\zeta_1$\\
   where $\zeta_0, \zeta_1, \zeta_2, \zeta_3$ is the dual basis of $S_0, S_1, S_2, S_3.$\\
   These elements form a basis for $A^!$ and in particular $dimA^!_m=(^4_m).$
   We also have the following Poincaré duality: $(A^!_m)^\ast\cong A^!_{4-m}.$
\end{propo}
If we consider $K_m(A)$ as free $A$-module, the Koszul complex of the algebra $A$ has the following form:
\begin{equation}\label{kcs}
0\longrightarrow A\stackrel{\bullet \textsf{t}}{\longrightarrow}S^4\stackrel{\bullet N}{\longrightarrow}S^6\stackrel{\bullet M}{\longrightarrow}S^4\stackrel{\bullet \textsf{x}}{\longrightarrow}S\stackrel{\varepsilon}{\longrightarrow}\mathbb C\longrightarrow 0.
\end{equation}
where: \\
$\bullet \textsf{x}$ is a right multiplication by $\textsf{x}=(S_0, S_1, S_2, S_3)^T;$\\
$\bullet M$ is the right multiplication by a matrix $M$ obtained from the relations $f_{0i}, f_{jk}$ of $A:$\\
$$M=\left(
  \begin{array}{cccc}
    -S_1 & S_0 & -\alpha_1S_3 & -\alpha_1S_2 \\
    S_1 & S_0 & S_3 & -S_2 \\
    -S_2 & -\alpha_3S_3 & S_0 & -\alpha_2S_1 \\
    S_2 & -S_3 & S_0 & S_1 \\
    -S_3 & -\alpha_3S_2 & -\alpha_3S_1 & S_0 \\
    S_3 & S_2 & -S_1 & S_0 \\
  \end{array}
\right);$$\\
$\bullet N$ is the right multiplication by the matrix $N$ given by:
$$N=\frac{1}{2}\left(
  \begin{array}{cccccc}
    2S_1 & 2S_2 & 2S_3 & 0 & 0 & 0 \\
    0 & (1-\alpha_2)S_3 & -(1+\alpha_3)S_2 & 2S_0 & (1+\alpha_2)S_3 & -(1-\alpha_3)S_2 \\
    -(1+\alpha_1)S_3 & 0 & (1-\alpha_3)S_1 & -(1-\alpha_1)S_3 & 2S_0 & (1+\alpha_3)S_1 \\
    (1-\alpha_1)S_2 & -(1+\alpha_2)S_1 & 0 & (1+\alpha_1)S_2 & -(1-\alpha_2)S_1 & 2S_0 \\
  \end{array}
\right);$$\\
and $\bullet \textsf{t}$ is the right multiplication by $t=(S_0, S_1, S_2, S_3).$
We denote by $\Delta$ the element $1\in K_4(A)$ and by $\Pi$ the element $(S_0, S_1, S_2, S_3)\in K_3(A).$ We have $b(\Delta)= b(\Pi)=0$ and:
$$q(\Pi)=3(S_2\otimes S_3\otimes f_{01}+S_3\otimes S_1\otimes f_{02}+S_1\otimes S_2\otimes f_{03}+S_0\otimes S_1\otimes f_{23})\in A^{\otimes 4};$$
$$q(\Delta)=1\otimes q(\Pi)\in A^{\otimes 5}.$$
Hence $\Pi$ and $\Delta$ define respectively elements of $HH_3(A)$ and $HH_4(A)$ which we will also denote using the same letters.
\section{Hochschild homology of the Sklyanin algebra $A$}
Assume that $\alpha_1, \alpha_2, \alpha_3$ are "generic". Formally, they generate the field $L=\mathbb Q(\alpha_1, \alpha_2, \alpha_3;\alpha_1+\alpha_2+\alpha_3+\alpha_1\alpha_2\alpha_3=0).$\\
We have:
$$
\begin{array}{lcl}
A&=& \mathbb C(\alpha_1, \alpha_2, \alpha_3)\langle S_0, S_1, S_2, S_3\rangle/(I_2)\\
&\cong&\mathbb C\otimes_{L}(L\langle S_0, S_1, S_2, S_3\rangle/(J_2)),
\end{array}
$$
where $(J_2)$ is the two-sided ideal of $L\langle S_0, S_1, S_2, S_3\rangle$ generated by $f_{0i}, \ f_{jk}.$
The we have the following commutative diagram:
{\small
$$\xymatrix{A^{\otimes n}\displaystyle{\ar[d]^b}\displaystyle{\ar[r]}&
\mathbb C\otimes_{L}\left(L\langle S_0, S_1, S_2, S_3\rangle/(J_2)\right)^{\otimes n}\displaystyle{\ar[d]^{1_{\mathbb C}\otimes_{L}b}}
\\
A^{\otimes (n-1)}\displaystyle{\ar[r]}&\mathbb C\otimes_{L}\left(L\langle S_0, S_1, S_2, S_3\rangle/(J_2)\right)^{\otimes (n-1)}}$$
}
where $b$ is the Hochschild boundary.\\
Therefore
$$HH_\star(A)\cong\mathbb C\otimes_{L}HH_\star(L\langle S_0, S_1, S_2, S_3\rangle/(J_2)),$$
as graded $\mathbb C$-vector spaces.\\
The morphism $L\longrightarrow\mathbb Q(\beta_1, \beta_2)((h)),$ $\alpha_1\mapsto \beta_1h^2,$ $\alpha_2\mapsto \beta_2h^2,$ $\alpha_3\mapsto =-(\beta_1h^2+\beta_2h^2)/(1+\beta_1h^2\beta_2h^2)=\beta_3h^2+O(h^3),$ with $\beta_3:=-\beta_1-\beta_2,$ defines an injection.\\
Let us set $k_0=\mathbb{Q}(\beta_1, \beta_2).$
\begin{propo}
The algebra $k_0((h))\langle S_0, S_1, S_2, S_3\rangle/(I)$ is isomorphic to $k_0((h))\otimes_{L}(L\langle S_0, S_1, S_2, S_3\rangle/(J_2))$ as graded $k_0((h))$-algebras,
while the vector space $$HH_\star(k_0((h))\langle S_0, S_1, S_2, S_3\rangle/(I))$$ is isomorphic to $$k_0((h))\otimes_{L}HH_\star(L\langle S_0, S_1, S_2, S_3\rangle/(J_2))$$
as graded $k_0((h))$-vector spaces.\\
Here $(I)$ is the two-sided ideal of $k_0((h))\langle S_0, S_1, S_2, S_3\rangle/(I)$ generated by:
$$\begin{array}{l}
F_{01}=[S_0, S_1]-\beta_1h^2(S_2S_3+S_3S_2);\\
F_{02}=[S_0, S_2]-\beta_2h^2(S_3S_1+S_1S_3);\\
F_{03}=[S_0, S_3]-(\beta_3h^2+O(h^3))(S_1S_2+S_2S_1);\\
F_{jk}=[S_j, S_k]-(S_0S_i+S_iS_0),
\end{array}$$
and $(i,j,k)$ is a cyclic permutation of $(1,2,3).$
\end{propo}
In particular we have\\
 {\small $HH_\star(L\langle S_0, S_1, S_2, S_3\rangle/(I_2))=L\otimes_{k_0((h))}HH_\star(k_0((h))\langle S_0, S_1, S_2, S_3\rangle/(I)).$} \\
 Therefore the computation of Hochschild's homology of algebra $A$ is equivalent to finding the homology of algebra $k_0((h))\langle S_0, S_1, S_2, S_3\rangle/(I)).$\\
 Let us introduce new variables by setting $x_0=-h^{-1}S_0;$ and $x_i=S_i$ for $i=1, 2, 3.$ We denote by $A_h$ the $k_0((h))$-algebra generated by $x_0, x_1, x_2, x_3$ with the relations
 $$\begin{array}{l}
g_{01}=[x_0, x_1]+\beta_1h(x_2x_3+x_3x_2);\\
g_{02}=[x_0, x_2]+\beta_2h(x_3x_1+x_1x_3);\\
g_{03}=[x_0, x_3]+(\beta_3h+O(h^2))(x_1x_2+x_2x_1);\\
g_{jk}=[x_j, x_k]+h(x_0x_i+x_ix_0),
\end{array}$$
where $(i,j,k)$ is a cyclic permutation of $(1,2,3).$
We also denote by $(I)$ the two-sided ideal of $k_0((h))\langle x_0, x_1, x_2, x_3\rangle$ generated by $g_{0i},  g_{jk}.$
$$A_h=k_0((h))\langle x_0, x_1, x_2, x_3\rangle/(I).$$
\subsection{Filtration on $A_h$}
Set $\mathcal A_h=k_0[[h]]\langle x_0, x_1, x_2, x_3\rangle/(I).$\\
$A_h=k_o((h))\otimes_{k_0[[h]]}\mathcal A_h$ has the same Poincaré series as $\mathcal R=k_0[x_0, x_1, x_2, x_3]=k_o\otimes_{k_0[[h]]}\mathcal A_h.$ Then $\mathcal A_h$ is a flat $k_0[[h]]$-module and the map $\mathcal A_h\longrightarrow A_h=k_o((h))\otimes_{k_0[[h]]}\mathcal A_h, a\mapsto 1\otimes a$ is an injection morphism. For this morphism, the image of $h$ is $h.$ Since $h$ is inversible in $A_h,$ $h$ is a nonzero divisor in $\mathcal A_h.$\\
The $h$-adic filtration $F$ on $\mathcal A_h$ can be extended to a filtration $F$ on $A_h$ such that the associated graded ring $gr_F(A_h)=k_0[x_0, x_1, x_2, x_3][h, h^{-1}].$ The filtred algebra $A_h$ is not complete. However each homogeneous component $(A_h)_n$, $n\in\mathbb N$, is complete.
 \subsection{The algebra $A_h$ seen as a deformation}
 Let $p$ be the projection  $p : \mathcal A_h\longrightarrow\mathcal R=k_0\otimes_{k_0[[h]]}\mathcal A_h, a\mapsto 1\otimes a.$ $Kerp=h\mathcal A_h=F_1\mathcal A_h.$\\
 Consider the commutator bracket $[\cdot, \cdot] : \mathcal A_h\times \mathcal A_h\longrightarrow\mathcal A_h.$ The terms $[x_i, x_j]$ can be written as formal series in $h$ with no constant coefficients in $\mathcal A_h.$\\ Therefore the map $h^{-1}[\cdot, \cdot] : \mathcal A_h\times \mathcal A_h\longrightarrow\mathcal A_h$ is a biderivation which satisfies the Jacobian identity.\\
 Let us consider the bracket: $$\{\cdot, \cdot\} : \mathcal R\times\mathcal R\longrightarrow\mathcal R$$
  defined by
  \begin{equation}
  \begin{array}{l}
  \{x_0, x_i\}=-2\beta_ix_jx_k \\ \\
  \{x_j, x_k\}=-2x_0x_i
  \end{array}
  \end{equation}
   where $(i, j, k)$ is a cyclic permutation of $(1, 2, 3).$\\
This is nothing but the Sklyanin Poisson bracket which it is the Jacobian Poisson structure given by:
\begin{equation}\label{e1}
P_1=x_1^2+x_2^2+x_3^2
\end{equation}
\begin{equation}\label{f1}
P_2=x_0^2+J_1x_1^2+J_2x_2^2+J_3x_3^2
\end{equation}
where $\beta_i=J_j-J_k$ and $(i, j, k)$ is a cyclic permutation of $(1, 2, 3).$\\
 We have the following commutative diagram:
 $$\xymatrix{\mathcal A_h\times \mathcal A_h\displaystyle{\ar[r]^{\ \ \frac{1}{h}[\cdot, \cdot]}}\displaystyle{\ar[d]^{p\times p}}&
\mathcal A_h\displaystyle{\ar[d]^{p}}\\
\mathcal R\times\mathcal R\displaystyle{\ar[r]^{\ \ \{\cdot, \cdot\}}}&\mathcal R}$$
\subsection{Hochschild homology of algebra $A_h$}
The filtration $F$ on $A_h$ can be extended to a filtration $F$ on $C(A_h)$ such that the associated grated ring $gr_F(C(A_h))=k_0[x_0, x_1, x_2, x_3][h, h^{-1}].$ \\
Let us consider the spectral sequence $E^r_{p, q}$ associated to the filtration $F$ on $C(A_h).$\\
Let $A_{p, q}^\infty=Ker\{d : F_pC_{p+q}\longrightarrow F_pC_{p+q-1}\}$ and consider the usual map of the spectral theory:
$$\Phi^r_{p, q} : A_{p, q}^\infty\longrightarrow E_{p, q}^r.$$
\begin{thm}
The spectral sequence $E^r$ associated with the filtration $F$ converges to the Hochschild homology $HH_\star(A_h)$ of $A_h.$
\end{thm}
\begin{proof}This is a direct consequence of the fact that for all $n\in\mathbb N$, $(C(A_h)_n, b)$ (the sub-complex of $(C(A_h), b)$ which is formed by the homogeneous components of degree $n$ for the initial graduation (weight graduation) $A_h$) is a complex offinite dimensional $k_0$-vector spaces and the filtration $F$ on $(C(A_h)_n$ is complete
\end{proof}
Since the graded ring $gr_F(A_h)$ associated with the filtration $F$ on $A_h$ is a polynomial algebra with coefficients in the ring $k_0[h, h^{-1}]$, by the Hochschild-Kostant-Rosenberg theorem, we have for all $n\in \mathbb N$, a quasi-isomorphism of $k_0[h, h^{-1}]$-modules:
$$\begin{array}{lll}
C_n(gr_F(A_h))&\longrightarrow&\Omega^n_{gr_F(A_h)_{|k_0[h, h^{-1}]}}\\
r_0\otimes r_1\otimes\cdots\otimes r_n&\mapsto&\frac{1}{n!}r_0dr_1\wedge\cdots\wedge dr_n.\\
\end{array}$$
Here $\Omega^n_{gr_F(A_h)_{|k_0[h, h^{-1}]}}$ is the $k_0[h, h^{-1}]$-module of differential forms of degree $n$ of $gr_F(A_h)$ on $k_0[h, h^{-1}],$ with zero differential. Hence $HH_ngr_F(A_h))\cong\Omega^n_{gr_F(A_h)_{|k_0[h, h^{-1}]}}.$ Under this quasi-isomorphism, the Connes's coboundary corresponds to the de Rham differential. On the other hand, we have a canonical isomorphism of $k_0[h, h^{-1}]$-modules $\Omega^\bullet_{gr_F(A_h)_{|k_0[h, h^{-1}]}}\cong \Omega^\bullet_{\mathcal R_{|k_0}}\otimes_{k_0}k_0[h, h^{-1}].$ We will denote $\Omega^\bullet_{\mathcal R_{|k_0}}$ by $\Omega^\bullet(\mathcal R).$\\
We borrow the following commutative diagram from Brylinski's paper (~\cite{bry}):
$$\xymatrix{E^1_{-n}=HH_n(gr_F(A_h))\displaystyle{\ar[r]^{\ \ \cong\ \ \  }}\displaystyle{\ar[d]^{d^1}}&\Omega^n(\mathcal R)\otimes_{k_0}k_0[h, h^{-1}]\displaystyle{\ar[d]^{\partial\otimes\cdot h}}\\
E^1_{-n+1}=HH_{n-1}(gr_F(A_h))\displaystyle{\ar[r]^{\ \ \ \cong \ \   }}&\Omega^{n-1}(\mathcal R)\otimes_{k_0}k_0[h, h^{-1}]}$$
where $d^1$ is the differential which computes the second term $E^2$ of the spectral sequence, $\cdot h$ is the multiplication by $h$ and $\partial$ is the boundary Poisson operator associated with the Sklyanin Poisson bracket:
$$\partial_n(F_0dF_1\wedge...\wedge dF_n)=\sum(-1)^{i+1}\{F_0, F_i\}dF_1\wedge...\wedge \widehat{dF_i}\wedge ...\wedge dF_n$$
$$+\sum (-1)^{i+j}F_0d\{F_i,F_j\}\wedge dF_1\wedge...\wedge \widehat{dF_i}\wedge...\widehat{dF_j}\wedge...\wedge dF_n$$\\
where $F_0,...,F_n\in\mathcal R.$\\
Using the isomorphism $E^1_{p, q}\cong\Omega^{p+q}(\mathcal R)\otimes_{k_0}h^{-p},$ the first term of this spectral sequence can be explained as follow:
{\tiny
$$\begin{array}{cccccccccccccc}
&\cdot & \cdot & &\cdot & &\cdot & &\cdot & &\cdot & &\cdot\\
&\cdot & \cdot & &\cdot & &\cdot & &\cdot & &\cdot & &\cdot\\
&\cdot & \cdot & &\cdot & &\cdot & &\cdot & &\cdot & &\cdot\\
\cdots&\Omega^{1}(\mathcal R)\otimes h^{-3}&\longrightarrow\mathcal R\otimes h^{-2}&\longrightarrow &0&\cdots &0&\cdots&0&\cdots&0&\cdots&q=-2\\ \\
\cdots&\Omega^{2}(\mathcal R)\otimes h^{-3}&\longrightarrow\Omega^{1}(\mathcal R)\otimes h^{-2}&\longrightarrow &\mathcal R\otimes h^{-1}&\cdots &0&\cdots&0&\cdots&0&\cdots&q=-1\\ \\
\cdots&\Omega^{3}(\mathcal R)\otimes h^{-3}&\longrightarrow\Omega^{2}(\mathcal R)\otimes h^{-2}&\longrightarrow &\Omega^{1}(\mathcal R)\otimes h^{-1}&\cdots &\mathcal R\otimes 1&\cdots&0&\cdots&0&\cdots&q=0\\ \\
\cdots&\Omega^{4}(\mathcal R)\otimes h^{-3}&\longrightarrow\Omega^{3}(\mathcal R)\otimes h^{-2}&\longrightarrow &\Omega^{2}(\mathcal R)\otimes h^{-1}&\cdots &\Omega^{1}(\mathcal R)\otimes 1&\cdots&\mathcal R\otimes h&\cdots&0&\cdots&q=1\\ \\
\cdots&0&\longrightarrow\Omega^{4}(\mathcal R)\otimes h^{-2}&\longrightarrow &\Omega^{3}(\mathcal R)\otimes h^{-1}&\cdots &\Omega^{2}(\mathcal R)\otimes 1&\cdots&\Omega^{1}(\mathcal R)\otimes h&\cdots&\mathcal R\otimes h^2&\cdots&q=2\\
&\cdot & \cdot & &\cdot & &\cdot & &\cdot & &\cdot & &\\
&\cdot & \cdot & &\cdot & &\cdot & &\cdot & &\cdot & &\\
&\cdot & \cdot & &\cdot & &\cdot & &\cdot & &\cdot & &\\
&p=3 & p=2& &p=1& & p=0&& p=-1&&p=-2&&\\
\end{array}$$}
The second term $E^2$ of the spectral sequence is given by the homology of the lines with respect to the differential $\partial\otimes \cdot h.$ Since the multiplication by $h$ is a $k_0[h, h^{-1}]$-isomorphism, to have this second term, we only have to find the Poisson homology:
\begin{equation}\label{phc}
0\longrightarrow\Omega^4(\mathcal R)\stackrel{\partial_4}{\longrightarrow}\Omega^3(\mathcal R)\stackrel{\partial_3}{\longrightarrow}\Omega^2(\mathcal R)\stackrel{\partial_2}{\longrightarrow}\Omega^1(\mathcal R)\stackrel{\partial_1}{\longrightarrow}\mathcal R
\end{equation}
The proposition \ref{sph} gives this Poisson homology in the generic case.
\begin{propo}
The spectral sequence associated with the filtration $F$ degenerates at $E^2$.
\end{propo}
\begin{proof}
This spectral sequence degenerates at $E^2$ if the map $\Phi_{p, q}^2$ is surjective for all $p, q\in\mathbb Z.$ But the columns of $E^2$ are the same up to a multiplication by $h.$ Thus we only have to give a proof for $\Phi_{0,q}^2, \ q\in\mathbb Z.$ \\
\begin{itemize}
  \item $E_{0,0}^2\cong PH_0(\mathcal R, \partial)$ is the quotient of $\mathcal R.$ Let $v\in E^2_{0,0}$. Then $v$ can be lifted to an element u of $\mathcal R\cong F_0A_h/F_{-1}A_h.$ On other hand, $u$ can be lifted to an element $U$
of $F_0A_h$ and $\Phi_{0,0}^2(U)=v.$\\
Then let $\widetilde{P}_i\in F_0(A_h), i=1, 2$ be an element which lifts  $P_i,$ where $P_1, P_2$ are Casimirs which give the Sklyanin Poisson bracket. Since $k_0[P_1, P_2]$ is the center of the Sklyanin Poisson algebra $(\mathcal R, \{\cdot, \cdot\})$, $k_0[\widetilde{P}_1, \widetilde{P}_2]$ is the center of algebra $A_h.$ We endow $A_h$ with a natural structure of $k_0[\widetilde{P}_1, \widetilde{P}_2]$-module.
  \item
Then using the result (\ref{sph}), as a $k_0[P_1, P_2]$-module,  $E_{0,1}^2\cong PH_1(\mathcal R, \partial)$ is generated by the class an element $f(P_1, P_2)d\psi\in\Omega^1(\mathcal R)$, where $d$ is the de Rham differential, $\psi\in\mathcal R.$  Let $\Psi\in F_0(A_h)$ be an element which lifts $\psi.$ $B(\Psi)\in F_0(A_h^{\otimes 2})$ and $b(B(\Psi))=-B(b(\Psi))=0.$ We have $\Phi_{0,1}^2(B(\Psi))=\overline{d \psi}.$
\item Let $v\in E^2_{0,2}\cong PH_2(\mathcal R, \partial).$ From (\ref{sph}), $v$ is the class of an element $f(P_1, P_2)dP_1\wedge d\psi,$ $f(P_1, P_2)\in k_0[P_1, P_2].$ Since $\overline{P_1d\psi}\in PH_1(\mathcal R, \partial)$, there exists $\Psi\in A^\infty_{0,1}$ such that $\Phi_{0,1}^2(\Psi)=\overline{P_1d\psi}.$ We have $b(\Psi)=0$ and therefore $b(B(\Psi))=0.$ $\Phi_{0,2}^2(B(\Psi))=\overline{dP_1\wedge d\psi}.$\\ Then $f(\widetilde{P}_1, \widetilde{P}_2)B(\Psi)$ lifts $v.$
\item Since the image of the Hochschild's cycle $\Pi\in A_h^{\otimes 4}$ in $gr_F(A_h^{\otimes 4})$ is the generator $\pi$ of $PH_3(\mathcal R, \partial)=E^2_{0,3}$, as a free $k_0[P_1, P_2]$-module, $f(\widetilde{P}_1, \widetilde{P}_2)\Pi$ lifts the element $f(P_1, P_2)\pi\in E^2_{0,3}.$
\item Similarly, since the image of the Hochschild's cycle $\Delta\in A_h^{\otimes 5}$ in $gr_F(A_h^{\otimes 5})$ is the generator $\delta$ of $PH_4(\mathcal R, \partial)=E^2_{0,4}$, as a free $k_0[P_1, P_2]$-module, $f(\widetilde{P}_1, \widetilde{P}_2)\Pi$ lifts the element $f(P_1, P_2)\delta\in E^2_{0,4}.$
    \end{itemize}
\end{proof}
\newpage
\begin{thm} In the generic case, the Hochschild Homology of $A_h$ is described as follow:
\begin{itemize}
  \item $HH_4(A_h)$ is a free $k_0((h))[\widetilde{P}_1, \widetilde{P}_2]$-module of rank $1$ generated by the homogeneous element $\Delta$ of degree $4.$
  \item $HH_3(A_h)$ is a free $k_0((h))[\widetilde{P}_1, \widetilde{P}_2]$-module of rank $1$ generated by the homogeneous element $\Pi$ of degree $4.$
  \item $HH_2(A_h)$ is a free $k_0((h))[\widetilde{P}_1, \widetilde{P}_2]$-module of rank $6$ generated by homogeneous elements of respective degrees $3, 3, 3, 3, 4, 4.$
  \item $HH_1(A_h)$ is a free $k_0((h))[\widetilde{P}_1, \widetilde{P}_2]$-module of rank $13$ generated by homogeneous elements of respective degrees $1, 1, 1, 1, 2, 2, 2, 2, 3, 3, 3, 3, 4.$
  \item $HH_0(A_h)$ is a free $k_0((h))[\widetilde{P}_1, \widetilde{P}_2]$-module of rank $7$ generated by homogeneous elements of respective degrees $0, 1, 1, 1, 1, 2, 2.$
\end{itemize}
\end{thm}
We can also deduce the following result:
\begin{coro}
As $k_0((h))$-vector spaces, the homological groups $HH_i(A_h)$ have the following Poincaré series:
$$\begin{array}{lcl}
P(P(HH_0(A_h),t)&=&\frac{2t^2+4t+1}{(1-t^2)^2} ;\\
&&\\
P(HH_1(A_h),t)&=&\frac{t^4+4t^3+4t^2+4t}{(1-t^2)^2};\\
&&\\
P(HH_2(A_h),t)&=&\frac{2t^4+4t^3}{(1-t^2)^2}; \\
&&\\
P(HH_3(A_h),t)&=&\frac{t^4}{(1-t^2)^2} ;\\
&&\\
P(HH_4(A_h),t)&=&\frac{t^4}{(1-t^2)^2}. \\
\end{array}$$
\end{coro}

\noindent
Laboratoire Angevin de Recherche en Math\'ematiques Universit\'e D'Angers D\'epartement de Math\'ematiques\\
\it{E-mail address} : \verb"pelap@math.univ-angers.fr "
\end{document}